%% file: main.tex
\documentclass[review,12pt]{elsarticle}

\usepackage{graphicx}
\usepackage{balance}
\usepackage{multirow}
\usepackage{placeins}
\usepackage{amssymb}
\usepackage{amsmath}
\usepackage{xcolor}

\usepackage{url}
\usepackage{graphicx}
\usepackage{indentfirst}
\usepackage{quoting}
\usepackage{booktabs}
\usepackage{subcaption}
\usepackage{caption}

\usepackage{tikz}
\tikzset{>=latex}

\usepackage{algpseudocode}
\usepackage[vlined,linesnumbered,ruled]{algorithm2e}
\SetKwProg{Procedure}{}{}{}
\SetKwProg{Function}{}{}{}

\usepackage{etoolbox}
\makeatletter
\patchcmd{\pprintMaketitle}
{\hrule\vskip12pt}
{\clearpage\hrule\vskip12\p@}
{}{}
\patchcmd{\pprintMaketitle}
{\hrule\vskip12pt}
{\hrule\clearpage}
{}{}
\makeatother

\journal{Applied Mathematical Modelling}

\begin{document}

\begin{frontmatter}

\title{Mono and Multi-objective Models for the Emergency Medical Service in the City of Belo Horizonte, Brazil}

\author[label1]{Mariana Mendes Guimarães}
\ead{mmguimaraesbr@gmail.com}
\author[label2]{Flávio Vinícius Cruzeiro Martins*}
\ead{flaviocruzeiro@cefetmg.br}
\author[label2]{Rodrigo Ferreira da Silva}
\ead{rfsilva.dcc.ufmg@gmail.com}
\address[label1]{Postgraduate Program in Mathematical and Computational Modeling, Centro Federal de Educação Tecnológica de Minas Gerais.}
\address[label2]{Computer Department, Centro Federal de Educação Tecnológica de Minas Gerais, Av. Amazonas 7675, Belo Horizonte, MG, Brazil, 30510-000}

\begin{abstract}
This work presents the facility location problem for maximum coverage in an emergency medical service called SAMU-BH, located in Belo Horizonte, Minas Gerais, Brazil. The purpose of the proposed model is to find the best locations for its ambulances in order to provide faster in-site assistance whenever necessary. First, it develops a mono-objective model to maximize the coverage rate, defined as the number of calls attended before a target time. Based on the Facility Location and Equipment Emplacement Technique (FLEET) model, a novel mathematical model is proposed, called FLEET-IC, which can work with independent coverages. Secondly, it proposes a multi-objective model to increase total coverage while reducing the number of installed bases. The results show that the mono-objective model's application increases the coverage rate up to 48\%, considering the busy fraction factor and up to 26\% without it. It also analyzes the effect on the coverage rate of the inclusion of additional ambulances.
\end{abstract}

\begin{keyword}
Facility Location Problem, Integer Linear Programming, Mono and Multi-objective problems.
\end{keyword}

\end{frontmatter}

\section{Introduction}

The Facility Location Problem (FLP) is a classic problem in Operational Research that aims to find the best facilities locations to supply local demand. The choice of the best location is driven by an objective function and a set of constraints. Typical objective functions encompass minimizing operating costs, minimizing average travel time, or maximizing demand points coverage. Applications of FLP are frequently constrained on capacity, maximum reach time, or staff availability \citep{daskin:08, revelle:05}.

This work focus on the scenario of SAMU-BH, an Emergency Medical Service (EMS) of the city of Belo Horizonte, Minas Gerais, Brazil. The emergency service of SAMU-BH begins with a call to 192 received by a 24/7 Regulatory Center. As a call is received, a regulatory physician evaluates the urgency of the situation and gives the user initial instructions. If necessary, he calls an ambulance that provides on-site assistance. According to the victim's location, the nearest suitable ambulance is allocated to the service. As soon as the ambulance arrives, the staff initiate the in-site procedures and then decide whether or not to transport the victim to the nearest suitable hospital.

Each ambulance has a fixed point, called base, where it is conducted when not in service. According to the Brazilian Federal Law, each base requires a parking location for the ambulance, bathrooms, and rooms for the ambulance's staff. The better bases distribution throughout the territory, the greater the chances of an ambulance arriving in a reasonable time \citep{hale:03}. Therefore, each ambulance base position's best choice is directly related to the service's efficiency and better public resources management.

The response time is the main criterion for measuring the service's performance. It is determined as the time elapsed from receipt incident on the Regulatory Center until the team arrives at the victim's location. A call is considered covered if the response time complies with the health management area's limits set. Therefore, coverage is defined as the percentage of covered calls during a specific time. SAMU-BH offers two types of ambulances, Advanced Support Units (USAs), each with an intensive care unit, and Basic Support Units (USBs), structured to ensure first-aid procedures. SAMU-BH management considers calls covered when attended in ten minutes or less for USAs and eight minutes or less for USBs.

\begin{figure}[!htb]
    \centering
    \includegraphics[width=0.5\textwidth]{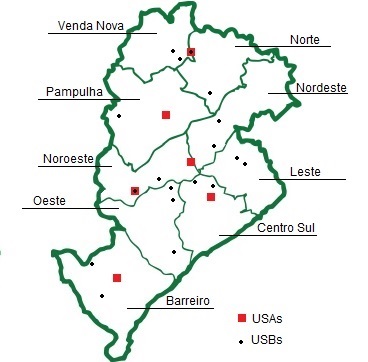}
    \caption{Approximate location of the SAMU-BH Bases in 2017. The division on the map corresponds to the administrative regions provided by the city of Belo Horizonte.}
    \label{fig:basesatuais1}
\end{figure}

From September 2015 to August 2016, SAMU-BH answered an average of 6,800 calls per month. In 2017, the service had 6 USAs and 21 USBs. All the 27 ambulances available were located in 22 bases (Figure \ref{fig:basesatuais1}); each base could receive one or more units of the same or different types. 

This work proposes a mono-objective model that aims to find the best location for SAMU-BH bases considering 1,527 distinct candidate points. It claims that if the models' application increases the coverage, the service of SAMU-BH is consequently improved. Based on the model Facility Location and Equipment Emplacement Technique (FLEET) \citep {schilling79}, a mathematical model adapted to the reality of SAMU-BH is developed. The main change is in the concept of coverage, which is broken down into advanced coverage (USAs) and basic coverage (USBs). For this reason, the model receives the name of FLEET-IC (Facility Location and Equipment Emplacement Technique with Independent Coverages). Additionally, a multi-objective model is proposed to reduce the number of bases and reduce installation and maintenance costs while increasing total coverage. Finally, new scenarios are also introduced to reproduce the impact by increasing the number of ambulances.

The remainder of this paper is structured as follows. The next section presents an overview of related works published in the literature. Section \ref{sec:methodology} shows how medical records are collected and used to generate new benchmark data sets. Next, in Section \ref{sec:mono}, it proposes an Integer Linear Programming model to optimize the ambulances' position and their bases' location. In Section \ref{sec:multi}, the proposed model is expanded to a multi-objective strategy that balances coverage of demands and the number of bases installed. Section \ref{sec:results} presents computational results. Finally, concluding remarks are presented in Section \ref{sec:conclusion}.

\section{Related Work}

~\citet{owen98} propose three categories for the strategies for facility location: the median problem, the central problem, and the coverage problem. The last one is the strategy used for the location of ambulances of SAMU-BH. The coverage problem is a typical strategy that applies to critical problems like emergency services location (e.g. fire brigade, ambulances). The main goal is to provide a solution to reach all users at a maximum acceptable time or distance, defined as the critical value. Some recent studies for SAMU can be found by ~\citep{junior11}, ~\citep{andrade14}, ~\citep{marques16} and ~\citep{guimaraes:18}.

The first mathematical model for the coverage problem is proposed by ~\citet{toregas71}. The essence of the coverage problem in facilities location is establishing a critical value that will be a reference to set if demand or user is covered.

Coverage problems can be divided into two distinct categories: set covering and maximum coverage problems. The first one, all demand, has to be covered, i.e., met within the established critical value. Also, in the set covering the number of facilities installed has to be minimized. However, 100\% coverage of demand becomes infeasible in many real-world applications, mainly due to budget constraints. The Maximal Covering Location Problem (MCLP) is an alternative in such cases that the goal is to maximize the covered demand within a pre-established critical value, considering the number of resources available. The first model for the maximum coverage problem is introduced by ~\citet{church74}.

The classic maximum coverage model's evolution is the FLEET (Facility Location and Equipment Emplacement Technique) model, presented by ~\citet{schilling79}. It considers different types of vehicles, and there is no hierarchy between them, i.e., they can be allocated independently. Besides, there is a distinction between base location points and ambulance allocation, which is essential for the later development of dynamic models. Due to these characteristics, the FLEET model is used as a reference for developing the mathematical models in this work.

A variant of the model FLEET is also introduced in \cite{RODRIGUEZ2020106522}. The Facility Location and Equipment Emplacement Technique model with Expected Covering (FLEET-EXC) model is used for fire stations location optimization. The authors present a mixed-integer linear programming model that considers vehicles' average utilization to compute expected demand coverage and an iterative procedure as a solving method. A Hypercube Queueing Model is used to compute the utilization of the vehicles. They present a case study on Concepcion province, Chile. They provide an analysis of the effects of locating and relocating fire stations for strategic decision making.

The facility location problem has become the preferred approach for dealing with emergency humanitarian logistical problems. Due to the increased number of natural and human-made disasters, \citet{BOONMEE2017485} propose a survey to show different works that use an exact algorithm, heuristic algorithm, or combined version of them as the primary approach to solving this problem. The survey examines the four main problems highlighted in the literature review: deterministic facility location problems, dynamic facility location problems, stochastic facility location problems, and robust facility location problems.

Due to external factors, the response time between any pair of vertices and vertices' demands in real-world facility location problems are usually indeterminate. \citet{ZHANG2017429} employ the uncertainty theory to address emergency service facilities' location problem under uncertainty.  First, the paper studies the uncertainty distribution of the covered demand associated with the covering constraint confidence level $\alpha$. After, the authors model the maximal covering location problem in an uncertain environment using different modeling ideas.

\citet{ANDERSSON2020103975} present three strategic cases for EMS planners as potential solutions to achieve high performance in this kind of service. In the first case, they measure the impact of replacing a closing down a local emergency room (ER) to adding an ambulance station and an ambulance in the area affected by the ER's closing. In the second case, they explore the potential of more effective utilization of the ambulances redirecting non-urgent calls for other transportation types. In the third case, they compared the average daily demand by taking time-varying demand into account. They extended the Maximum Expected Performance Location Problem with Heterogeneous Regions (MEPLP-HR) model to include multiple time periods.

Developing Countries have more difficulty to maintain a high level of healthcare. With this focus, \citet{doi:10.1287/opre.2019.1969} propose a robust optimization approach. It is used to optimize both the location and routing of emergency response vehicles, accounting for uncertainty in travel times and spatial demand characteristics of low- and middle-income countries. Throw a combination of a prediction-optimization framework with a simulation model and real data to provide an in-depth investigation; they can provide useful service analyses. They show the possibility of reducing the median average response time by roughly 10\%–18\% over the entire week and 24\%–35\% during rush hour.

\citet{10.1145/2001858.2001907} show an evolutionary metaheuristic for optimizing EMS applied to a real-world case in Argentina. The authors propose a bi-objective formulation, reducing service delay time and minimizing third-party medical vehicle use. The approach's results show maximizing the use of the available installed capacity, improving response time rates, and using a smaller number of resources.

A bi-objective robust program to design a cost-responsiveness efficient EMS system under uncertainty is presented by \citet{ZHANG20141033}. Because of uncertain parameters standard in EMS systems, they develop a robust optimization approach that determines the location of EMS stations and, at the same time, the number of EMS vehicles at each station. It is used to guarantee the balance between cost and responsiveness.

\citet{Wang2016} address a new two-stage optimization method for emergency supplies allocation problem with multi-supplier, multi-affected area, multi-relief, and multi-vehicle. The first stage is selecting a set of candidate relief suppliers and the number of various vehicles. The second stage deal with the decisions in terms of vehicle routing and relief allocation in every disaster scenario. It proposes a multi-objective optimization with goals to minimize the proportion of demand unsatisfied and response time of emergency reliefs and the whole process's total cost. Fuzzy random variables are used to deal with some unknown events.

To deal with emergency medical resources allocation and optimization scheduling of resources, \citet{Wen:2017:2156-7018:393} present a mathematical optimization model and an improved multi-objective algorithm based on an artificial bee colony algorithm. The authors report that the results show that the method achieves better time performance and better results than traditional multi-objective optimization algorithms. 

In a large scale emergency scenario, \citet{PAUL2017147} formulate a multi-objective hierarchical extension of the maximal covering location problem that seeks to maximize coverage of the population within a rapid response window while minimizing modifications to the existing structure. The authors analyze the effectiveness of the current and optimal locations of a set of existing regional assets maintained by the Department of Defense of the United States of America to respond to large-scale emergencies. Applying the $\epsilon$-constraint method, they show the trade-off between maximizing coverage and minimizing cost. The results show that it is possible to improve coverage by more than 15\%, and, with less than a 30\% modification of the entire structure, the coverage can exceed 98\%.

\citet{Hesam2020} also use a robust optimization approach to deal with the location-routing problem (LRP) under uncertainty for providing EMS during disasters. It proposes a bi-objective model to minimize relief time and the total cost, including location costs and route coverage costs (ambulances and helicopters). The authors propose a novel robust mixed-integer linear programming model to formulate the problem as an LRP. To solve the problem, it suggests a shuffled frog leaping algorithm (SFLA). The performance is compared using both the $\epsilon$-constraint method and the NSGA-II algorithm. The results indicate the efficiency of the SFLA within a proper computation time compared to the CPLEX solver as an exact method. They show that the SFLA has a better performance than the NSGA-II when compared in 4 quality indicators.

\section{Data Sets}
\label{sec:methodology}

 The collected data set of medical records comprises 29,048 occurrences, 12,777 for USAs, and 16,271 for USBs, registered between May 2016 and April 2017 for which an ambulance has been dispatched. Because of the time limitations, the sample comprising USBs occurrences for one-third of working days. SAMU-BH provided from paper all data used in this work, and now all of that is available in digital format\footnote{The approval to use such data for the matter of this research was granted by the Research Ethics Committee of the Municipal Health Department of Belo Horizonte (SMSA-BH) on November 1st, 2017, reference number 2,361,551.}. For each occurrence, there are the following attributes: a unique identifier, neighborhood and city of occurrence, ambulance commitment time, ambulance arrival and departure time from the place of occurrence, destination, the ambulance arrival time for hospital care, and rescue unit's release time.

A new data set is created based on the collected occurrences of SAMU-BH to test different scenarios. Each new data set consists of a set of occurrences generated by using \textit{Poisson} distribution considering the average of occurrences by the day of the week, time range, type of occurrence, type of ambulance; and, neighborhood where the call came from. 

\begin{figure}[ht!]
\begin{subfigure}[t]{0.48\textwidth}
  \centering
  \includegraphics[clip, trim=11cm 8.5cm 10cm 10cm, width=0.94\textwidth]{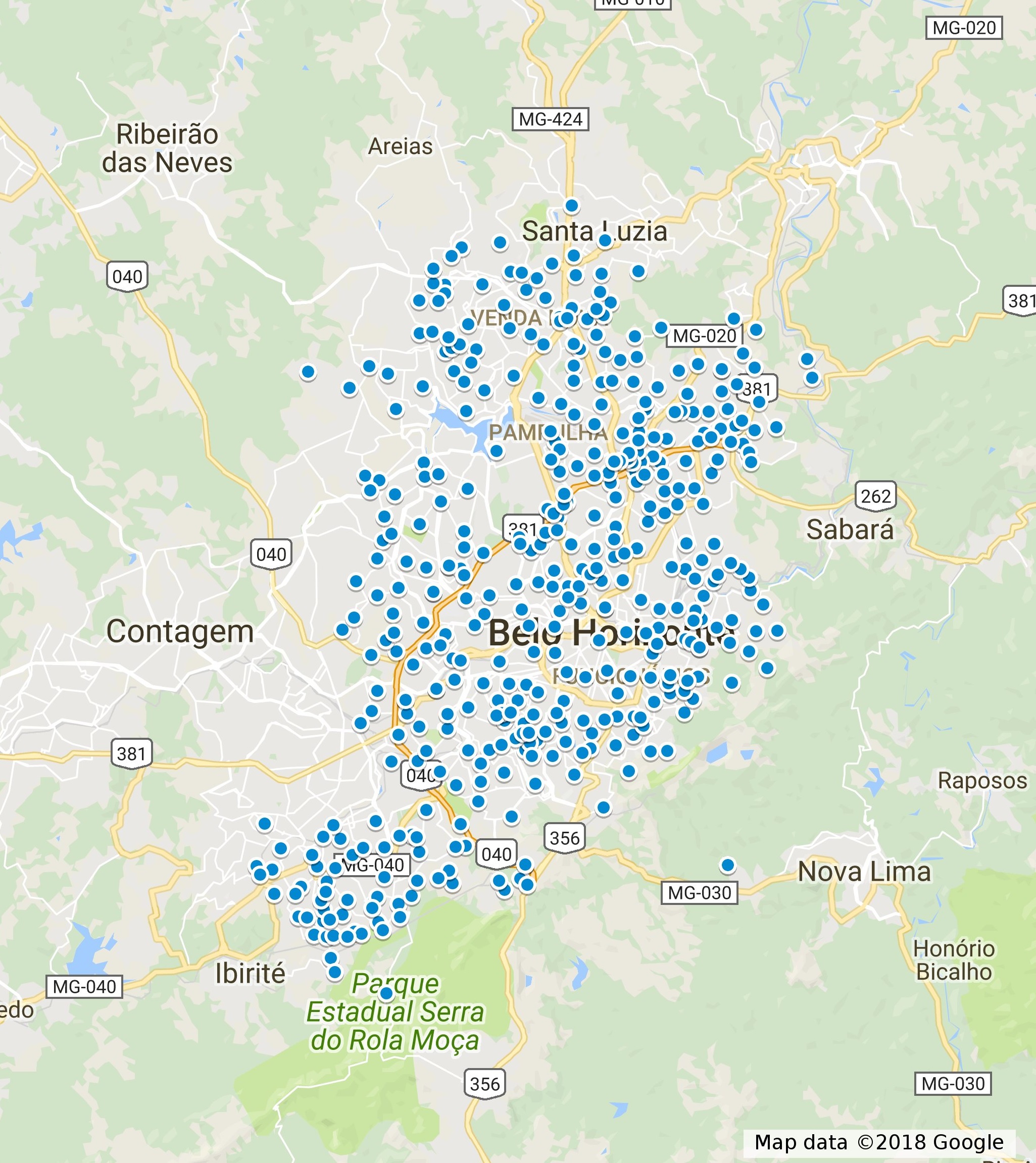}  
  \caption{Belo Horizonte neighborhoods mapped as 427 demand nodes. Each demand point is the centroid of the corresponding Belo Horizonte neighborhood.}
  \label{fig:demandPoints}
\end{subfigure}
\hfill
\begin{subfigure}[t]{0.48\textwidth}
  \centering
  \includegraphics[clip, trim=0.4cm 0cm 1.1cm 0cm, width=\textwidth]{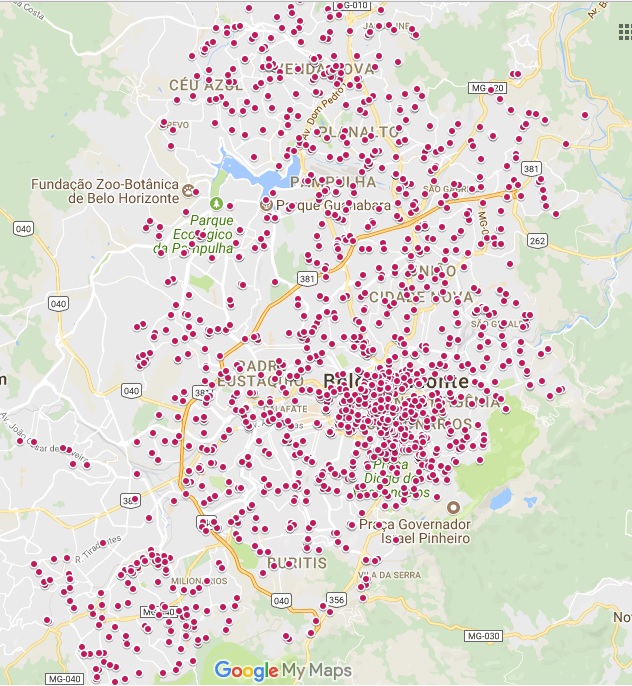}  
  \caption{1,527 candidates nodes for SAMU-BH base location in Belo Horizonte.}
  \label{fig:candidatesPoint}
\end{subfigure}
\caption{Demand and candidates nodes showed in the Belo Horizonte map using \textit{Google Maps}, 2017}
\label{fig:fig}
\end{figure}

As suggested by \citet{goldberg04}, all the demand locations are aggregated in the centroid of the corresponding neighborhood, resulting in 427 distinct points (Figure: \ref{fig:demandPoints}). This way, the service's response time is calculated considering the time to travel from the base where the dispatched ambulance is located and the occurrence neighborhood's centroid.

A set of 1,527 new points are mapped among hospitals, state buildings, parking lots, malls, and schools as candidates to new SAMU-BH base locations (Figure: \ref{fig:candidatesPoint}). According to the Brazilian Ministry of Health protocols, all mapped points have the structure to accommodate the facility installation ~\citep{brasil2014}.

The next section describes the mathematical model employed to obtain new optimized locations for the bases and ambulances.

\section{Mono-objective Model}
\label{sec:mono}

This section shows a new Integer Linear Programming model for SAMU-BH based on the FLEET \citep{schilling79} model. The model's main change is the concept of coverage split in coverage for USAs and coverage for USBs. For this reason, it names this new model FLEET-IC (Facility Location and Equipment Emplacement Technique with Independent Coverages). 

FLEET-IC has nine parameters, namely:

\begin{enumerate}
\item[$I$] set of demand nodes;
\item[$J$] set of candidates for facilities locations;
\item[$U$] set of rescue team types;
\item[$P^u$] the number of rescue teams available of type $u \in U$;
\item[$Q$] the number of bases to be installed;
\item[$S^u$] the limit time to attend a call for each rescue team type $u \in U$;
\item[$C_j$] the maximum amount of rescue teams for each location $j \in J$, also considering that each base can receive only one rescue teams of each type $u$;
\item[$d_{ji}$] the time to travel from a base installed in $j \in J$ to the demand node $i \in I$; and,
\item[$q_{iu}$] the average number of calls originating from each point $i \in I$ for each type of rescue team $u \in U$, based on historical data.
\end{enumerate}

FLEET-IC has only binary variables, that are described below:

\begin{enumerate}
\item[$z_j$] equal to 1 if a base is installed at location $j \in J$ and 0, otherwise;
\item[$x_j^u$] equal to 1 if a rescue team of type $u \in U$ is installed  for at location $j \in J$ and 0, otherwise; and,
\item[$y_{ji}^u$] equal to 1 if the location $i \in I$ is covered by a base installed at location $j \in J$ with a rescue team of type $u \in U$, such that $d_{ij} \leq S^u$;
\item[$k_{iu}$] equal to 1 if the location $i \in I$ is covered by rescue team of type $u \in U$ and 0, otherwise.
\end{enumerate}

\begin{equation}
\label{eq:modelopropriofo}
maximize \enspace coverage = \sum_{i \in I} \sum_{u \in U} q_{iu} k_{iu}
\end{equation}

\noindent subject to:

\begin{equation}
\label{eq:modeloproprior1}
\sum_{j \in J} x_j^u \leq P^u, \qquad \forall u \in U,
\end{equation}

\begin{equation}
\label{eq:modeloproprior2}
\sum_{u \in U} x_j^u \leq C_j, \qquad \forall j \in J,
\end{equation}

\begin{equation}
\label{eq:modeloproprior3}
\sum_{j \in J} z_j \leq Q,
\end{equation}

\begin{equation}
\label{eq:modeloproprior4}
\sum_{u \in U} x_j^u \leq U z_j,  \qquad \forall j \in J,
\end{equation}

\begin{equation}
\label{eq:modeloproprior5}
y_{ji}^u d_{ji} \leq S^u  x_j^u,  \qquad \forall i \in I, \forall j \in J, \forall u \in U, 
\end{equation}

\begin{equation}
\label{eq:modeloproprior6}
k_{iu} \leq \sum_{j \in J} y_{ji}^u,  \qquad \forall i \in I, \forall u \in U,
\end{equation}

\begin{equation}
\label{eq:modeloproprior7}
z_j, x_j^u, y_{ji}^u, k_{iu} \in \{0,1\}, \qquad \forall i \in I, \forall j \in J, \forall u \in U.
\end{equation}

The objective function (\ref{eq:modelopropriofo}) maximizes the number of covered calls. Historical data of calls on each demand point is considered for prioritizing locations with higher demand, and it is represented by $q_{iu}$.

Inequations (\ref{eq:modeloproprior1}), (\ref{eq:modeloproprior2}) and (\ref{eq:modeloproprior3}) define restrictions on capacity of amount of rescue teams available, amount of rescue team on each base and the total amount of bases that can be installed, respectively. Restriction (\ref{eq:modeloproprior4}) forbids a rescue team to be installed in a location without base installed. Restriction (\ref{eq:modeloproprior5}) matches the rescue team type with demand points needs, also considering the time limit to cover the call for the specific type of rescue team needed. A call is considered covered if there is at least one rescue team covering the call location, as imposed by restriction (\ref{eq:modeloproprior6}). Finally, (\ref{eq:modeloproprior7}) defines variables domains.

\section{Multi-objective Model}
\label{sec:multi}

Multi-objective optimization is a good alternative to handle problems with more than one objective. The implemented multi-objective approach works on two conflicting objectives to optimize the rescue teams' locations. While the number of covered calls is maximized, the number of bases installed is minimized. Both functions are normalized. Equations \ref{eq:mof1} and \ref{eq:mof2} show the function to maximize the coverage, $TxCob$, and to minimize the number of bases, $TxRedBases$,
\begin{equation}
\label{eq:mof1}
f_1(.) = TxCob = \frac{\sum_{i \in I} \sum_{u \in U} q_{iu} k_{iu}}{\sum_{i \in I} \sum_{u \in U} q_{iu}},
\end{equation}

\begin{equation}
\label{eq:mof2}
f_2(.) = TxRedBases = \frac{Q-\sum_{j \in J} z_{j}}{Q}.
\end{equation}

A linear scalarization is done to deal with the multi-objective formulation in the mathematical solver. This approach that uses the weighted sum method is called \textit{a posteriori}. Equation \ref{eq:mofo} shows the resulting objective function, 
\begin{equation}
Max \enspace z = \lambda(TxCob) + (1-\lambda)(TxRedBases).
\label{eq:mofo}
\end{equation}

The other restrictions remain the same as the mono-objective approach. Generating 101 distinct values for $\lambda$ ranging from 0 to 1, by incrementing 0.01. In this way, the optimization phase results in a wide range of solutions that compose the Pareto set, which can be analyzed and taken into account in the decision-making process.

\section{Busy Fraction}
\label{sec:busyfraction}

This section introduces the concept of the busy fraction of each advanced and basic rescue team unit. Busy fraction, represented by $q$, is estimated as described by \citet{daskin:82} via Equation \ref{eq:malp1rho1},

\begin{equation}
\label{eq:malp1rho1}
q = \frac{\overline{t} \sum_{i \in I} a_i}{24\sum_{j \in J} x_j} = \frac{\overline{t} \sum_{i \in I} a_i}{24P}.
\end{equation}


\input{txocupacaoA.tex}
\input{txocupacaoB.tex}

\noindent Which $\overline{t}$ is the average time in hours needed to attend a call, including response time fully, time on site, hospital drop-off, and time to return to base. The number of calls for each position is represented by $a_i$. Binary variable $x_j$ is set to true, if and only if there is a rescue team on position $j$. $ P$ represents the number of available rescue teams.

Tables \ref{tab:txocupacaoA} and \ref{tab:txocupacaoB} shows the busy fraction for USBs and USAs, respectively, with data ranging from May/2016 to April/2017. By using the average busy fraction found for each type of rescue team, it is possible to estimate the minimum number of rescue teams needed, $b$, which covers all demand points with a confidence level of $\theta$, calculated by Equation \ref{eq:malp1b1} as \citep{revelle:89},
\begin{equation}
\label{eq:malp1b1}
b = \left\lceil \frac{log(1-\theta)}{log (q)}\right\rceil.
\end{equation}
Table \ref{tab:calculob} shows values for $b$, for USAs, as $b_1$, and USBs, as $b_2$, with confidence levels of 0.95, 0.90 and 0.85.

\input{calculob.tex}

\section{Results}
\label{sec:results}

\begin{figure}[!htb]
    \centering
    \includegraphics[width=\textwidth]{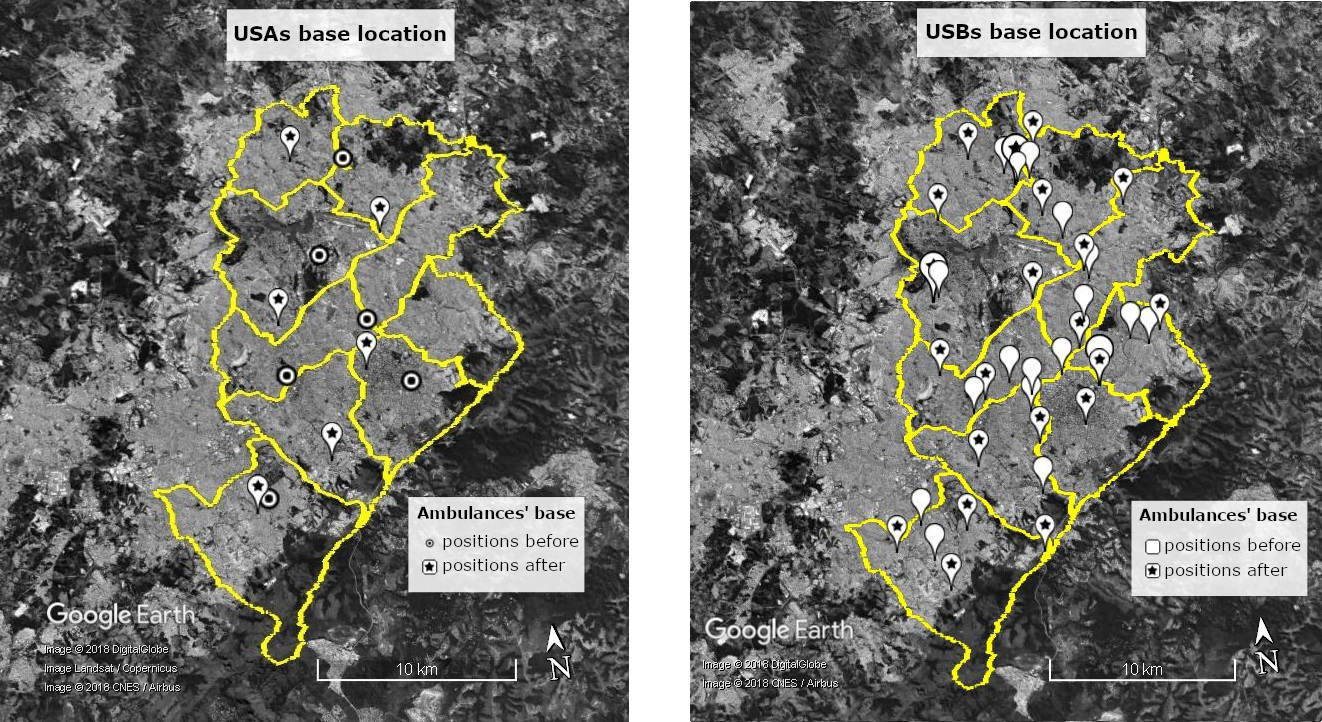}
    \caption{Current ambulances' base location versus Scenario 3 optimization.}
    \label{fig:result}
\end{figure}

The FLEET-IC model is run in IBM ILOG CPLEX$^{\textregistered}$ 12.7.1, using its default settings, on a computer with an Intel Core i5 5200U 2.2GHz processor with 8GB of RAM, on Windows. Six scenarios are considered. Table \ref{tab:resultado} present the results for the first three scenarios, with 26.736 calls, 6 USAs, and 21 USBs available. The first scenario simulates the current bases and rescue teams' positions of SAMU-BH represented in collected data. In the second scenario, the FLEET-IC model is used to optimize rescue teams' locations, keeping bases unchanged. Furthermore, in the third scenario, the model optimizes location bases and rescue teams using 1.527 possible points to install bases. Results show that the model increases coverage of calls in scenarios 2 and 3, by a margin of 10\% and 26\%, respectively. It is interesting to note that the SAMU in scenario 2 does not require any investment to implement the solution. It is possible because this scenario does not consider any new basis. Only reallocations can be made. The locations of the bases for each ambulance type before and after optimization are shown in Figure \ref{fig:result}. The division on the map corresponds to the administrative regions provided by the city of Belo Horizonte.
\input{tabelaresultados}

The fourth scenario considers multiple calls occurring at the same time. In this case, even if a base covers a call, the rescue team will attend to the demand only if it is available. The busy fraction presented in Section \ref{sec:busyfraction} is used to check whether the rescue team is available or not. The presented confidence levels of 0.95, 0.90, and 0.85. Results are shown in Table \ref{tab:resultado2} for 26,736 calls, 06 USAs and 21 USBs. When the FLEET-IC model is used, the coverage increases between 31\% and 48\%, taken into account all 1.527 possible base locations.

\input{tabelaresultados2}

The fifth and sixth scenarios are built by simulating the impact of introducing 1, 2, and 3 rescue teams on scenarios three and four. First, only the new rescue teams' location is optimized by keeping the remaining rescue teams in their original location. With extras USB and extras USA, the coverage is increased from 69\% to 77\%, 82\%, and 87\%, with one, two, or three ambulances of each type, respectively. However, suppose the location of original and new ambulances is optimized. In that case, the coverage increases from 69\% to 92\%, 94\%, and 96\%, with one, two, or three ambulances of each type, respectively.

\begin{figure}[!htb]
    \centering
    \includegraphics[width=1.0\textwidth]{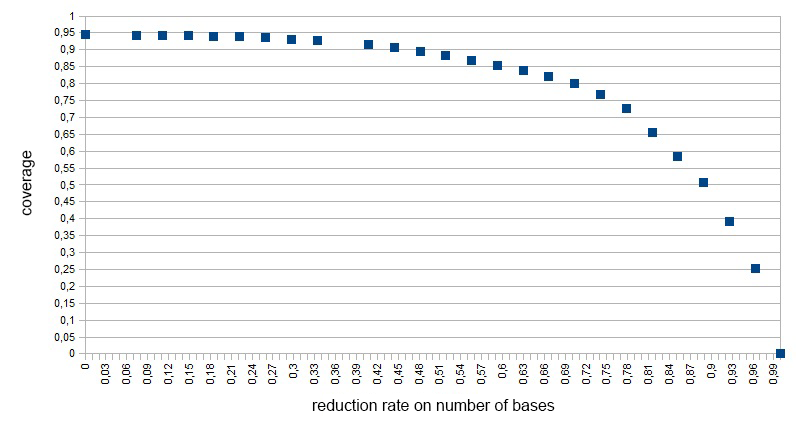}
        \caption{Pareto Set using 101 weight variations in Equation \ref{eq:mofo}.}
    \label{fig:grafpareto}
\end{figure}

The multi-objective optimization is an excellent tool for decision-makers in the public administration to improve the decision process. All tests for the multi-objective approach are performed according to scenario 3. The 101 weight variations, as previously described, generates a set of Pareto solutions with 26 non-dominated solutions, as shown in Figure \ref{fig:grafpareto}. 

As a decision-maker, SAMU-BH management can choose a solution accordingly. Among the criteria that can be analyzed, the following stand out:
i) variation in the percentage of coverage of advanced services; 
ii) variation in the percentage of coverage of essential services; 
iii) variation in the percentage of total coverage; iv) direct and indirect costs of installing a base; 
v) number of points covered; 
vi) allocation of more than one same type unit on the same basis according to the most significant demand nodes queue.

\begin{figure}[!htb]
    \centering
    \includegraphics[width=\textwidth]{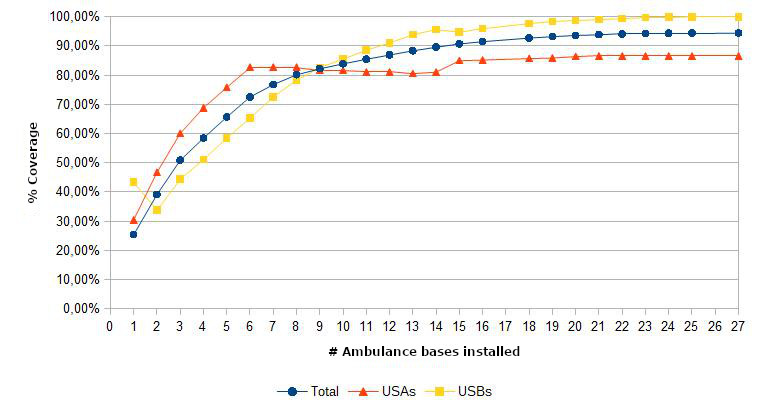}
    \caption{Coverage rates variation versus the number of ambulance bases installed.}
    \label{fig:covxnumamb}
\end{figure}

Figure \ref{fig:covxnumamb} presents a graph that shows the relationship between the call coverage rate and the number of installed bases.  It is important to note that the coverage model objective is to maximize the total number of calls covered, that is, the sum of calls made by USAs and USBs within response times. Therefore, in some cases, the model may prioritize one ambulance type to maximize the objective function's value.

The USAs coverage rate is not sensitive to the reduction of up to six bases. In other words, for 27 or 21 installed bases, the coverage of calls would remain the same for this kind of ambulance. As noted before, the coverage of some ambulance types can decrease in some situations.  It happens for 14 installed bases of USBs, which has covered a more significant number of calls (95.62\%) than the optimal installation of 15 bases (94.76\%). As the objective is to maximize the total coverage, the bases' optimal location points follow the criterion of privileging the nodes with the highest volume of calls, whether advanced or basic units.

When analyzing the total coverage, a reduction of 50\%  in the number of bases decreases the covered calls rate by 5\%. It can happen since some more distant neighborhoods generated a tiny number of calls. Consequently, it is often necessary to install a base at a new point to cover one or a few occurrences. For example, when $\lambda = 0.99$, there is a reduction of two bases. The coverage is decreased by only 0.04\% (11 calls are no longer covered).

\section{Concluding Remarks}
\label{sec:conclusion}

This work applies the Facility Location Problem for Maximum Coverage in an Emergency Medical Service called SAMU-BH. In this type of service, resources are limited and must be used efficiently. Thus, it proposes two models to optimize bases and rescue teams' locations, increase coverage, and support demand. 

Using the first proposed model, called FLEET-IC, the coverage is increased up to 48\% considering the busy fraction factor, and between 10\% and 26\%  without it. The results show the possibility to increase the coverage up to 10\% without any investment, using just ambulances relocation strategy. Considering additional rescue teams' acquisition, the coverage rate increases between 8\% and 27\%, depending on the parameters used.

Secondly, it proposes a multi-objective approach that maximizes demand coverage while minimizing the number of bases installed. Several optimal non-dominated solutions are obtained that balance coverage and number of bases. This information can support managers in the decision-making process to improve the service for the coming years.

Computational results certify Operational Research as a powerful tool to support health services management decisions and open new perspectives for further studies in the subject. This work shows that theory and practice can be put together to solve a real-life problem and effectively contribute to providing a better health service.

\section*{Acknowledgments}

Authors thank SAMU-BH, PPSUS, CAPES, CNPq and FAPEMIG for supporting this research.

\bibliographystyle{elsarticle-harv}
\bibliography{main}

\end{document}

%% file: txocupacaoA.tex
\begin{table}[!htb]
    \footnotesize
    \centering
    \caption{Busy Fraction for USAs}
    \label{tab:txocupacaoA}
    \begin{tabular}{lrrr}
        \toprule
\multicolumn{1}{c}{\textbf{}} & \multicolumn{2}{c}{\textbf{Total time (hours)}} & \multicolumn{1}{c}{\textbf{Busy Fraction}} \\
\multicolumn{1}{c}{\textbf{Month}} & \multicolumn{1}{c}{In service} & \multicolumn{1}{c}{Available} &  \multicolumn{1}{c}{\textbf{$q$}} \\
        \midrule
May-2016	&	763.15	&	4464	&	0.17 \\
Jun-2016	&	806.42	&	4320	&	0.19 \\
Jul-2016	&	699.27	&	4464	&	0.16 \\
Aug-2016	&	651.72	&	4464	&	0.15 \\
Sep-2016	&	617.52	&	4320	&	0.14 \\
Oct-2016	&	600.33	&	4464	&	0.13 \\
Nov-2016	&	592.82	&	4320	&	0.14 \\
Dec-2016	&	691.82	&	4464	&	0.15 \\
Jan-2017	&	634.63	&	4464	&	0.14 \\
Feb-2017	&	534.12	&	4032	&	0.13 \\
Mar-2017	&	636.02	&	4464	&	0.14 \\
Apr-2017	&	618.80	&	4320	&	0.14 \\
        \midrule
&	(sum) 7846.60 &	(sum) 52560 &	(avg) 0.15  \\
        \bottomrule
    \end{tabular}
\end{table}

%% file: txocupacaoB.tex
\begin{table}[!htb]
    \footnotesize
    \centering
    \caption{Busy Fraction for USBs}
    \label{tab:txocupacaoB}
    \begin{tabular}{lrrr}
        \toprule
\multicolumn{1}{c}{\textbf{}} & \multicolumn{2}{c}{\textbf{Hours}} & \multicolumn{1}{c}{\textbf{Busy Fraction}} \\
\multicolumn{1}{c}{\textbf{Month}} & \multicolumn{1}{c}{In service} & \multicolumn{1}{c}{Available} &  \multicolumn{1}{c}{\textbf{$q$}} \\
        \midrule
May-2016	&	960.07	&	4032	&	0.24	\\
Jun-2016	&	857.08	&	3528	&	0.24	\\
Jul-2016	&	885.75	&	3528	&	0.25	\\
Aug-2016	&	948.15	&	4032	&	0.24	\\
Sep-2016	&	942.57	&	3528	&	0.27	\\
Oct-2016	&	856.07	&	3528	&	0.24	\\
Nov-2016	&	733.77	&	3528	&	0.21	\\
Dec-2016	&	903.07	&	4032	&	0.22	\\
Jan-2017	&	832.18	&	3528	&	0.24	\\
Feb-2017	&	816.40	&	3528	&	0.23	\\
Mar-2017	&	928.68	&	3528	&	0.26	\\
Apr-2017	&	912.03	&	4032	&	0.23	\\
        \midrule
&	(sum) 10575.82	&(sum) 44352	& (avg)	0.24	\\
        \bottomrule
    \end{tabular}
\end{table}

%% file: calculob.tex
\begin{table}[!htb]
    \footnotesize
    \centering
    \caption{Number of USAs and USBs to cover a location}
    \label{tab:calculob}
    \begin{tabular}{rr|rrr|rrr}
        \toprule
\multirow{2}{*}{\textbf{$\theta$}} & \multirow{2}{*}{$log(1-\theta)$} & \multicolumn{3}{|c}{\textbf{USA}} & \multicolumn{3}{|c}{\textbf{USB}} \\

\multicolumn{1}{c}{} & \multicolumn{1}{c}{} & \multicolumn{1}{|c}{$log(q_1)$} & \multicolumn{1}{c}{$\frac{log(1-\theta)}{log(q_1)}$} & \multicolumn{1}{c}{\textbf{$b_1$}} & \multicolumn{1}{|c}{$log(q_2)$} & \multicolumn{1}{c}{$\frac{log(1-\theta)}{log(q_2)}$} & \multicolumn{1}{c}{\textbf{$b_2$}} \\
       \midrule
\textbf{0.95}	&	-1.30	&	-0.82	&	1.58	&	\textbf{2}	&	-0.62	&	2.10	&	\textbf{3}	\\
\textbf{0.90}	&	-1.00	&	-0.82	&	1.21	&	\textbf{2}	&	-0.62	&	1.61	&	\textbf{2}	\\
\textbf{0.85}	&	-0.82	&	-0.82	&	1.00	&	\textbf{1}	&	-0.62	&	1.33	&	\textbf{2}	\\
\textbf{0.80}	&	-0.70	&	-0.82	&	0.85	&	\textbf{1}	&	-0.62	&	1.13	&	\textbf{2}	\\
        \bottomrule
   \end{tabular}
\end{table}

%% file: tabelaresultados.tex
\begin{table}[!htb]
    \footnotesize
    \centering
    \caption{Scenarios 1, 2 and 3}
    \label{tab:resultado}
    \begin{tabular}{lrrr}
        \toprule
\multicolumn{1}{c}{} & \multicolumn{1}{c}{Scenario 1} & \multicolumn{1}{c}{Scenario 2} & \multicolumn{1}{c}{Scenario 3} \\
        \midrule
Bases installed: & 22 & 22 & 27 \\
Computational time:	&	-	&	11''  &	31'' \\
Objective function:	&	18.521	&	21.015 &	25.242\\
Coverage for USAs:	&	61\%	&	74\% &	87\% \\
Coverage for USBs: 	&	75\%	&	82\% &	100\% \\
\hline
\textbf{Overall coverage}:	&	\textbf{68\%}	&	\textbf{78\%} &	\textbf{94\%}\\
      \bottomrule
    \end{tabular}
\end{table}

%% file: tabelaresultados2.tex
\begin{table}[!htb]
    \footnotesize
    \centering
    \caption{Scenarios 1 and 4 for different confidence levels for busy fraction}
    \label{tab:resultado2}
    \begin{tabular}{lrrrrrr}
        \toprule
\multicolumn{1}{c}{} & \multicolumn{3}{c}{Scenario 1} & \multicolumn{3}{c}{Scenario 4} \\
        \midrule
Confidence level $\theta$: 	&	0.95	&	0.9	&	0.85	&	0.95	&	0.90	&	0.85	\\
USAs to cover: 	&	2	&	2	&	1	&	2	&	2	&	1	\\
USBs to cover: 	&	3	&	2	&	2	&	3	&	2	&	2	\\
Bases installed: 	&	22	&	22	&	22	&	27	&	26	&	25	\\
Computational time: 	&	-	&	-	&	-	&	3703''	&	178''	&	98''	\\
Objective function:	&	3.802	&	8.390	&	14.416 &	16.567	&	19.706	&	22.819	\\
Coverage for USAs &	7\%	&	7\%	&	61\% &	58\%	&	58\%	&	87\%		\\
Coverage for USBs &	20\%	&	49\%	&	49\%	 &	64\%	&	85\%	&	85\%	\\
\hline
\textbf{Overall coverage:}	&	\textbf{14\%}	&	\textbf{31\%}	&	\textbf{54\%}	&	\textbf{62\%}	&	\textbf{74\%}	&	\textbf{85\%}\\
      \bottomrule
    \end{tabular}
\end{table}